\documentclass[11pt]{article}
\usepackage{amsmath}

\usepackage[margin=1in, centering]{geometry}

\usepackage{placeins}
\usepackage{verbatim}

\usepackage{authblk}
\usepackage{url}
\usepackage{amsthm}

\usepackage[hyphenbreaks]{breakurl}

\usepackage[dvips]{graphicx}
\usepackage{color}
\usepackage{amsfonts,amssymb,amsmath,}
\usepackage{mathrsfs}
\usepackage[T1]{fontenc}
\usepackage{ae}

\usepackage{bbm}

\newtheorem{theorem}{Theorem}

\newtheorem{lem}{Lemma}
\newtheorem{defn}[lem]{Defintion}

\numberwithin{equation}{section}
\numberwithin{table}{section}
\numberwithin{figure}{section}

\renewcommand{\(}{\left(}
\renewcommand{\)}{\right)}

\allowdisplaybreaks[4]

\begin{document}
\title{Large bias for integers with prime factors in arithmetic progressions}
\author{Xianchang Meng}
\date{}
\maketitle
\begin{abstract}
We prove an asymptotic formula for the number of integers  $\leq x$ which can be written as the product of $k ~(\geq 2)$ distinct primes $p_1\cdots p_k$ with each prime factor in an arithmetic progression $p_j\equiv a_j \bmod q$, $(a_j, q)=1$ $(q \geq 3, 1\leq j\leq k)$. For any $A>0$, our result is uniform for $2\leq k\leq A\log\log x$. Moreover, we show that, there are large biases toward certain arithmetic progressions $(a_1 \bmod q, \cdots, a_k \bmod q)$, and such biases have connections with Mertens' theorem and the least prime in arithmetic progressions. 
\end{abstract}

\let\thefootnote\relax
\let\thefootnote\relax

\section{Introduction and statement of results}

For any $k\geq 2$, $q\geq 3$, and integers $(a_j, q)=1 ~(1\leq j\leq k)$, we consider the number of integers $\leq x$ which can be written as product of $k$ distinct primes $p_1 p_2 \cdots p_k$ with $p_j\equiv a_j\bmod q$ $(1\leq j\leq k)$. Here when we count the number of such integers, we allow any ordering of the prime factors.

Ford and Sneed \cite{Ford-S} investigated subtle biases in the distribution of the product of two primes in different arithmetic progressions subject to the Generalized Riemann Hypothesis (GRH) and the Linear Independence conjecture (LI) on the imaginary parts of the zeros of Dirichlet L-functions.  Recently, the author \cite{Meng} generalized their results to study the
bias of numbers composed of $k$ ($\ge 2$) prime factors (either
counted with or without multiplicity) in two different arithmetic progressions. 
For each fixed $k$ and $q$, different arithmetic progressions contain
virtually the same number of such integers below $x$, indeed, under the GRH, the number
equals $\frac{1}{\phi(q)}\frac{x}{\log x}\frac{(\log\log x)^{k-1}}{(k-1)!} + O(x^{1/2+o(1)})$. 


It is reasonable to expect that these integers break up very evenly,
with errors of size $O(x^{1/2+o(1)})$, when
one specifies further which arithmetic progression modulo $q$ each prime factor
lies in.  However, this is not the case.  Dummit, Granville and Kisilevsky \cite{DGK}
showed that there is a very large bias for
the odd integers $p_1p_2 \le x$ with two prime factors
satisfying $p_1 \equiv p_2 \equiv 3\pmod{4}$.  More precisely, they showed that
$$ \frac{\#\{p_1p_2\leq x: p_1\equiv p_2\equiv 3\bmod{4} \}}{\frac{1}{4}\#\{p_1 p_2\leq x \}}=1+\frac{c+o(1)}{\log\log x},$$
for some positive constant $c$.  The authors exhibit a similar bias
for products of 2 primes, where $\chi_q(p_1)=\chi_q(p_2)=\eta$,
$\chi_q$ is a quadratic Dirichlet character with fixed conductor $q$,
and $\eta\in \{-1,1\}$.  If $q$ is allowed to grow with $x$, they
further conjecture that the bias may be a bit larger.  
 Recently, Hough \cite{Hou} confirmed their conjecture and showed that, for $\eta \in\{-1, 1 \}$, there exist many $q\leq x$ for which 
$$ \frac{\# \{p_1 p_2\leq x: \chi_q (p_1)=\chi_q(p_2)=\eta \}}{\frac{1}{4} \# \{p_1p_2\leq x: (p_1 p_2, q)=1 \}}~ \text{is at least as large as}~ 1+\frac{\log\log\log x+O(1)}{\log\log x}.$$
 On the other hand, Moree \cite{More} considered all the integers with every prime factor from the same arithmetic progression $a \bmod q$, and proved that there is a large bias towards certain residue classes $a \bmod q$. 


In this paper, we generalize the large bias results found in \cite{DGK} to products of any $k\geq 2$ primes and any fixed modulus $q\geq 3$, and prove uniform estimates in a large range of $k$. For any fixed $A>0$ and fixed $q\geq 3$, we prove an asymptotic formula uniformly for $2\leq k\leq A\log\log x$ for the number of integers $p_1\cdots p_k\leq x$ with $p_j\equiv a_j \bmod q ~(1\leq j\leq k)$ . We show that, there are large biases for some arithmetic progressions $(a_1 \bmod q, \ldots, a_k \bmod q)$,  and such phenomenon has connections with Mertens theorem and the least prime in arithmetic progressions. 

Let $\mathbf{a}:=(a_1, a_2, \dots, a_k)\in (\mathbb{Z}/q\mathbb{Z})^k$ with $(a_j, q)=1$ for all $ 1\leq j\leq k$. One may regard the
vector $\mathbf{a}$ as an unordered $k$-tuple, or as a multiset. Denote 
$$M_k(x; \mathbf{a}):=\#\{n\leq x: n=p_1 p_2\cdots p_k, p_j \text{~distinct primes}, p_j\equiv a_j \bmod q, (a_j, q)=1, 1\leq j\leq k \},$$
where the prime factors $p_j$ can be in any order, and 
$$S_k(x):=\#\{n\leq x: n=p_1 p_2 \cdots p_k, (p_1 p_2 \cdots p_k, q)=1, p_j ~\text{distinct primes}, 1\leq j\leq k \}.$$ 
Let $\chi$ be a Dirichlet character modulo $q$, and $\chi_0$ be the principal character modulo $q$. Denote
$$C(q, a):=\lim_{x\rightarrow\infty}\Bigg(\phi(q) \sum_{\substack{p\leq x\\ p\equiv a \bmod q}} \frac{1}{p} -\sum_{\substack{p\leq x\\ (p, q)=1}} \frac{1}{p}\Bigg)=\sum_{\chi\neq \chi_0}\bar{\chi}(a)\sum_{p}\frac{\chi(p)}{p}.$$
We will see in our theorems that this constant $C(q, a)$ reflects the bias in our problem. 

Our first result is for the special case when $a_1=a_2=\cdots =a_k=a$. In other words, all the $k$ prime factors are from the same residue class $a \bmod q$. 

\begin{theorem}\label{thm-same}
	Let $q\geq 3$ be fixed, and  $\mathbf{a}=(a, a, \dots, a)$, $(a, q)=1$. We have the following results. 
	
	i) For fixed $k\geq 2$, 
\begin{equation*}
\frac{M_k(x; \mathbf{a})}{\frac{1}{\phi^k(q)}S_k(x)}= 1+\frac{(k-1)C(q, a)}{\log\log x}+O_{q, k}\(\frac{1}{(\log\log x)^2} \).
\end{equation*}

ii) If $k=o(\log\log x)$, as $x\to\infty$, 
\begin{equation*}
\frac{M_k(x; \mathbf{a})}{\frac{1}{\phi^k(q)}S_k(x)}= 1+\frac{(k-1)(C(q, a)+o(1))}{\log\log x}.
\end{equation*}

iii) For fixed $A>0$ and $k\sim A\log\log x$,   we have, as $x\rightarrow\infty$,
\begin{equation*}
\frac{M_k(x; \mathbf{a})}{\frac{1}{\phi^k(q)}S_k(x)}\sim \prod_p \frac{1+\frac{A\phi(q)\boldsymbol{\mathbbm{1}}_{p\equiv a \bmod q}(p)}{p}}{1+\frac{A\chi_0(p)}{p}}.
\end{equation*}

\end{theorem}

\noindent\textbf{\textit{Remark 1.}} If $k$ is fixed, by Lemmas \ref{lem-Mk-same-pf} and \ref{lem-Sk-pf},  $M_k(x; \mathbf{a})$ and $\frac{1}{\phi^k(q)}S_k(x)$ have main terms of the same order which is $\frac{1}{\phi^k(q)}\frac{x}{\log x} \frac{(\log\log x)^{k-1}}{(k-1)!}$ with different secondary terms and hence the bias is determined by the constant $C(q, a)$. 
Thus we see that, as $k$ increases, the bias will become larger and larger.

For $k\sim A\log\log x$, the main terms of $M_k(x; \mathbf{a})$ and $\frac{1}{\phi^k(q)}S_k(x)$ have the same order of magnitude but with different coefficients. One may compare this with the result of Moree \cite{More} who showed that the counting function $N(x; q, a):=\#\{n\leq x: p|n\Rightarrow p\equiv a \bmod q \}$ satisfies $N(x; q, a)\sim B_{q, a} x/ (\log x)^{1-1/\phi(q)}$ for some positive constant $B_{q, a}$ depending on $q$ and $a$, and in particular, $N(x; 4, 3)\geq N(x; 4, 1)$ holds for all $x$.

\vspace{0.6em}


For the general case, assume there are $l$ distinct values $b_1, \dots, b_l$ in the coordinates of $\mathbf{a}$. Fix $l$, for each $1\leq j\leq l$, let $k_j$ be the number of prime factors congruent to $b_j \bmod q$. Then $\sum_{j=1}^l k_j=k$. 
\begin{theorem}\label{thm-general-quotient}
Let $q\geq 3$ be fixed. Then, for fixed $k\geq 2$, 
\begin{align*}
\frac{M_k(x; \mathbf{a})}{\frac{1}{\phi^k(q)}\frac{k!}{k_{1}!k_2!\cdots k_l!} S_k(x)}= 1+\frac{k-1}{\log\log x} \frac{1}{k}\sum_{j=1}^k C(q, a_j)+O_{q, k, l}\(\frac{1}{(\log\log x)^2} \).
\end{align*}
Moreover,  for  fixed $l$ and fixed $A>0$, assume $k=\sum_{j=1}^l k_j \sim A\log\log x$ and $e_j:=\lim_{x\rightarrow\infty}\frac{k_j}{\log\log x}$ exists for every $1\leq j\leq l$. Then as $x\rightarrow\infty$, 
\begin{equation}\label{thm-general-formu-uniform}
\frac{M_k(x; \mathbf{a})}{\frac{1}{\phi^k(q)}\frac{k!}{k_{1}!k_2!\cdots k_l!} S_k(x)}\sim \prod_p \frac{\prod_{j=1}^l \(1+\frac{\phi(q)e_j \boldsymbol{\mathbbm{1}}_{p\equiv b_j \bmod q}(p)}{p}\)}{1+\frac{A\chi_0(p)}{p}},
\end{equation}
where $\sum_{j=1}^{l}e_j=A$. 
\end{theorem}
\noindent{\textbf{\textit{Remark 2.}}} 
In the general case $\mathbf{a}=(a_1, \dots, a_k)$, there are $\frac{k!}{k_{1}!k_2!\cdots k_l!}$ orderings of the numbers $a_1, \ldots, a_k$. 

\noindent{\textbf{\textit{Remark 3.}}} For $k\sim A\log\log x$, if the coordinates of $\mathbf{a}$ cover all the reduced residue classes modulo $q$ and all the $e_j$'s are the same, then  the right side of (\ref{thm-general-formu-uniform}) is exactly $1$.

\subsection{Mertens theorem and the least prime in arithmetic progressions}
The constant $C(q, a)$, which affects the biases in our theorems, is related to the classical Mertens theorem (\cite{Hardy}, \S 22.8)
and the Mertens theorem \cite{Mert} for arithmetic progressions, that 
\begin{equation}\label{Mertens-thm}
\sum_{p\leq x}\frac{1}{p}=\log\log x+\gamma+B+O\(\frac{1}{\log x}\),
\end{equation}
and if $(a, q)=1$, 
\begin{equation}\label{Mertens-thm-AP}
\sum_{\substack{p\leq x\\p\equiv a\bmod q}}\frac{1}{p}=\frac{\log\log x}{\phi(q)}+M(q, a)+O\(\frac{1}{\log x}\), 
\end{equation}
where $\gamma$ is Euler's constant,  $B:=\sum_p\(\log\(1-\frac{1}{p}\)+\frac{1}{p}\)$ is Mertens' constant, and $M(q, a)$ is a number depending on $q$ and $a$.  Languasco and Zaccagnini \cite{Langu} investigated the value of $M(q, a)$ and other related constants. 
By (\ref{Mertens-thm}), (\ref{Mertens-thm-AP}), and the orthogonality of Dirichlet characters, letting $x\rightarrow\infty$, we get
\begin{equation}\label{C-bias-to-M-AP}
C(q, a)=\phi(q)M(q, a)-\gamma-B+\sum_{p|q}\frac{1}{p}, 
\end{equation}
\begin{equation}\label{Merten-C-sum}
\sum_{\substack{a \bmod q\\ (a, q)=1}}M(q, a)=\gamma+B-\sum_{p|q}\frac{1}{p}.
\end{equation}
 Hence the value of $M(q, a)$ determines how the bias behaves. 
 
 In particular,  with the values of $M(q, a)$ calculated by Languasco and Zaccagnini \cite{Langu-web}, by (\ref{C-bias-to-M-AP}), we have
 $$C(3, 2)\approx 0.641945, \quad C(3, 1)\approx -0.641945;$$
 $$C(4, 3) \approx 0.334981, \quad C(4, 1)\approx -0.334981;$$
 $$C(7, 2)\approx 1.83747, \quad C(7, 5)\approx 0.159006, \quad C(7, 6)\approx -0.946269;$$
 $$C(13, 3)\approx 2.68478, \quad C(13, 6)\approx -0.846522, \quad C(13, 8)\approx -1.31962.$$
Here the interesting phenomenon is that 2 is a quadratic residue modulo 7, while 5 and 6 are quadratic non-residues modulo 7; 3 is a quadratic residue modulo 13, while 6 and 8 are quadratic non-residues modulo 13. There is no consistent preference for either quadratic non-residue classes or quadratic residue classes modulo $q$. 

The above phenomenon is different from the biases among products of $k$ primes studied in \cite{Ford-S} and \cite{Meng}. 
Using the similar method as in \cite{Meng}, one can show that, under the GRH and LI, the integers $n=p_1\cdots p_k$, which are products of exactly $k$ distinct primes, have preference for either quadratic non-residues or quadratic residues, depending on the parity of $k$. 

The biases in Theorems \ref{thm-same} and \ref{thm-general-quotient} ultimately stem from the fact that $M(q, a)$ is heavily dependent on the least prime $p(q, a)$ in the arithmetic progression $a\bmod q$. Pomerance \cite{Pomera} and Norton \cite{Norton} independently showed that  
\begin{equation}\label{Pome}
\sum_{\substack{p\leq x\\p\equiv a \bmod q}}\frac{1}{p}-\frac{\log\log x}{\phi(q)}=\frac{1}{p(q, a)}+O\( \frac{\log 2q}{\phi(q)}\), 
\end{equation}
where the implied constant is uniform for all $q$, $a$, and $x\geq q$. 

In Theorem \ref{thm-general-quotient}, we allow any ordering of the primes $p_j$ $(1\leq j\leq k)$, and hence the constant $\frac{1}{k}\sum_{j=1}^k C(q, a_j)$ represents the bias. One may ask, for which $\mathbf{a}=(a_1, \cdots, a_k)$, this constant is 0? Trivially, by (\ref{C-bias-to-M-AP}) and (\ref{Merten-C-sum}), if $\mathbf{a}$ covers every element of the reduced residue class modulo $q$ the same number of times, $\frac{1}{k}\sum_{j=1}^k C(q, a_j)=0$. But we don't know if the converse is true. Alternatively, by (\ref{C-bias-to-M-AP}) and (\ref{Merten-C-sum}), we may consider the distribution of the values of $M(q, a_j)$ $(1\leq j\leq k)$. By (\ref{Pome}), it is reasonable to conjecture that all the $M(q, a_j)$'s are distinct and that, except in the trivial case, they are linearly independent over $\mathbb{Q}$. Hence, we propose the following open problem.

\smallskip

\noindent\textbf{Open Problem.} 
 Is the trivial case the only case for which $\frac{1}{k}\sum_{j=1}^k C(q, a_j)=0$?
 
 The answer is yes if the numbers $\sum_{p}\frac{\chi(p)}{p}$ are linearly independent over algebraic numbers or the numbers $\bar{\chi}(a)\sum_p \frac{\chi(p)}{p}$ are linearly independent over $\mathbb{Q}$ for all $\chi\neq \chi_0 \bmod q$ . These values are close to $\log L(1, \chi)$ or $\bar{\chi}(a)\log L(1, \chi)$. Baker, Birch, and Wirsing \cite{Baker} showed that if $(q, \phi(q))=1$ then the numbers $L(1, \chi)$ are linearly independent over $\mathbb{Q}$ for all non-principal characters $\bmod ~q$. Moreover, they proved that the numbers $L(1, \chi)$ are linearly independent over algebraic numbers for non-trivial even characters $\bmod ~q$ (see also  \cite{Murtys} Corollary 2 or  \cite{Murty-Rath} Corollary 25.6). For any odd Dirichlet character $\chi$,  it is known (\cite{Murty-Rath} Lemma 25.7) that the number $L(1, \chi)$ is an algebraic multiple of $\pi$.

\section{Lemmas and Preparations}

\begin{lem}[\cite{Kara}, Chapter IX, \S 2, Theorem 2, \cite{Dave}, page 96, (12)]\label{lem-zerofree}
The Dirichlet $L$-function	$L(s, \chi)$ has no zeros in the domain
	$$ \Re(s)=\sigma\geq 1-\frac{c_1}{\log q(|t|+2)}, $$
	for some constant $c_1>0$, except a possible simple real zero close to 1 when $\chi$ is real, which is called a Siegel zero. If $\chi$ is real, there exists an effective constant $c_2>0$ such that $L(\sigma, \chi)\neq 0$ in the range 
	$$\sigma>1-\frac{c_2}{\sqrt{q}\log^2 q}.$$
	
\end{lem}

We need the following terminologies (Part II. Chapter 5.2, \cite{Tene}).

\begin{defn}\label{defn-P}
Let $z\in \mathbb{C}$, $c_0>0$, $0<\delta\leq 1$, $M>0$. We say that a Dirichlet series $F(s)$ has the property $\mathcal{P}(z; c_0, \delta, M)$ if the Dirichlet series 
$G(s; z):= F(s)\zeta(s)^{-z}$
can be analytically continued to the region $\sigma\geq 1-c_0/(\log(|t|+2))$, and in this region, 
$|G(s; z)|\leq M (1+|t|)^{1-\delta}.$	
\end{defn}

\begin{defn}\label{defn-T}
	We say $F(s)$ has type $\mathcal{T}(z, w; c_0, \delta, M)$, 
	if $F(s)=\sum_{n\geq 1} a_n/n^s$ has property $\mathcal{P}(z; c_0, \delta, M)$, and there exists a sequence of non-negative real numbers $\{b_n \}_{n=1}^{\infty}$ such that $|a_n|\leq b_n$, and the series $\sum_{n\geq 1} b_n/n^s$ satisfies $\mathcal{P}(w; c_0, \delta, M)$ for some complex number $w$, 
\end{defn}

\begin{lem}[\cite{Tene}, Part II, Theorem 5.2]\label{lem-Dirichlet-z}
	Let $F(s):=\sum_{n\geq 1} a_n/n^s$ be a Dirichlet series of type $\mathcal{T}(z, w; c_0, \delta, M)$. For $x\geq 3$, $N\geq 0$, $A>0$, $|z|\leq A$, and $|w|\leq A$, we have
	\begin{equation*}
	\sum_{n\leq x} a_n=x(\log x)^{z-1}\left\{\sum_{0\leq n\leq N} \frac{u_n(z)}{(\log x)^n}+O(MR_n(x)) \right\}, 
	\end{equation*}
	where
	\begin{equation*}
	u_n(z):=\frac{1}{\Gamma(z-n)} \sum_{l+j=n} \frac{1}{l!j!} G^{(l)}(1; z) \gamma_j(z),
	\end{equation*}
	\begin{equation*}
	G^{l}(s; z):=\frac{\partial^l}{\partial s^l}G(s, z), \quad \gamma_j(z):=\frac{d^j}{ds^j}\( \frac{\{(s-1)\zeta(s)\}^z}{s}\), 
	\end{equation*}
	and 
	\begin{equation}\label{RN-defn}
	R_N(x)=e^{-c_1 \sqrt{\log x}}+\(\frac{c_2N+1}{\log x} \)^{N+1},
	\end{equation}
	for some constants $c_1$ and $c_2$ depending at most on $c_0$, $\delta$, and $A$. 
\end{lem}

\begin{lem}\label{lem-Tene-Ck}
	Let $a_z(n)$ be an arithmetic function depending on a complex parameter $z$ and $a_z(n)=\sum_{k=0}^{\infty} c_k(n)z^k$ in the disk $|z|\leq A$. Suppose there exists a function $h(z)$ holomorphic for $|z|\leq A$, and a quantity $R(x)$, independent of $z$, such that, for $x\geq 3$ and $|z|\leq A$, we have
	\begin{equation*}
	\sum_{n\leq x} a_z(n)=x(\log x)^{z-1}\left\{zh(z)+O_A(R(x)) \right\}.
	\end{equation*}
If $|h''(z)|\leq B_1$ for $|z|\leq A$, then uniformly for $x\geq 3$, $1\leq k\leq A\log\log x$, we have
\begin{equation}\label{for-lem-Tene-Ck}
C_k(x)=\frac{x}{\log x} \frac{(\log\log x)^{k-1}}{(k-1)!}\left\{h\(\frac{k-1}{\log x}\)+O_A\(\frac{B_1(k-1)}{(\log\log x)^2}+\frac{\log\log x}{k}R(x) \) \right\}.
\end{equation}
If we suppose $|h^{(4)}(z)|\leq B_2$ for $|z|\leq A$, then uniformly for $x\geq 3$, $3\leq k\leq A\log\log x$, we have
\begin{align}
C_k(x)=\frac{x}{\log x} \frac{(\log\log x)^{k-1}}{(k-1)!}&\bigg\{h(0)+\frac{k-1}{\log\log x} h'(0)+\frac{(k-1)(k-2)}{(\log\log x)^2}g\(\frac{k-3}{\log\log x} \) \nonumber \\
&+  O_A\(\frac{B_2 (k-1)(k-2)(k-3)}{(\log\log x)^4} +\frac{\log\log x}{k}R(x) \) \bigg\}, \label{for-lem-Tene-Ck-D2}
\end{align}
where $$g(z)=\int_0^1 h''(tz)(1-t)dt.$$
\end{lem}
\noindent \textit{\textbf{Proof.}}  Formula \eqref{for-lem-Tene-Ck} is a special case of Theorem 6.3 Part II in \cite{Tene}. 

For all $r\leq A$, the main term in \eqref{for-lem-Tene-Ck-D2} is from 
$$I:=\frac{x}{\log x}\frac{1}{2\pi i} \oint_{|z|=r} h(z)\frac{e^{z\log\log x}}{z^{k}}dz=\frac{x}{\log x}\frac{1}{2\pi i} \oint_{|z|=r} (h(0)+zh'(0)+z^2 g(z))\frac{e^{z\log\log x}}{z^{k}}dz, $$
where $g(z)=\int_0^1 h''(tz)(1-t)dt$. When $k\leq A\log\log x$, choose $r_j=\frac{k-j}{\log\log x}$ ($1\leq j\leq 3$), we see that 
\begin{align}
I=&\frac{x}{\log x} \frac{1}{2\pi i} \oint_{|z|=r_1} h(0)\frac{e^{z\log\log x}}{z^{k}}dz+\frac{x}{\log x} \frac{1}{2\pi i} \oint_{|z|=r_2} h'(0)\frac{e^{z\log\log x}}{z^{k-1}}dz\nonumber\\
&+\frac{x}{\log x} \frac{1}{2\pi i} \oint_{|z|=r_3} g(z)\frac{e^{z\log\log x}}{z^{k-2}}dz\nonumber\\
=&\frac{x}{\log x} \frac{(\log\log x)^{k-1}}{(k-1)!}\bigg\{ h(0)+\frac{k-1}{\log\log x}h'(0)\bigg\}+\frac{x}{\log x}\frac{1}{2\pi i}\oint_{|z|=r_3}g(z)\frac{e^{z\log\log x}}{z^{k-2}}dz.\label{pf-lem-Tene-main}
\end{align}
Next, we examine the last integral in \eqref{pf-lem-Tene-main}. Since we assume $|h^{(4)}(z)|\leq B_2$ for $|z|\leq A$, we have
\begin{align*}
g(z)=&g(r_3)+(z-r_3)g'(r_3)+(z-r_3)^2 \int_0^1 (1-t)g''(r_3+t(z-r_3)) dt\\
=&g(r_3)+(z-r_3)g'(r_3)+O\(B_2|z-r_3|^2 \).
\end{align*}
Thus, the last integral in \eqref{pf-lem-Tene-main} equals
\begin{align}
&\frac{x}{\log x}\Bigg\{ \frac{g(r_3)}{2\pi i}\oint_{|z|=r_3} \frac{e^{z\log\log x}}{z^{k-2}}dz+ \frac{1}{2\pi i}\oint_{|z|=r_3} (z-r_3) \frac{e^{z\log\log x}}{z^{k-2}}dz\nonumber\\
&\qquad+O\(B_2\int_0^{2\pi} |e^{i\alpha}-1|^2 e^{r_3\log\log x \cos\alpha} r_3^{5-k}d\alpha  \) \Bigg\}\nonumber\\
&=\frac{x}{\log x}\left\{g(r_3) \frac{(\log\log x)^{k-3}}{(k-3)!}+\frac{(\log\log x)^{k-4}}{(k-4)!}-r_3\frac{(\log\log x)^{k-3}}{(k-3)!} +O\( B_2\frac{(\log\log x)^{k-5}}{(k-4)!}\)\right\}\nonumber\\
&=\frac{x}{\log x} \frac{(\log\log x)^{k-3}}{(k-3)!}\left\{ g\(\frac{k-3}{\log\log x}\) +O\(\frac{B_2 (k-3)}{(\log\log x)^2} \) \right\}.\label{pf-lem-Tene-erro}
\end{align}
The error term $O\(R(x)\log\log x/k \)$ is the same as that in the proof of \eqref{for-lem-Tene-Ck}. Combing \eqref{pf-lem-Tene-main} and \eqref{pf-lem-Tene-erro}, we get the desired result. \qed

\bigskip

We need some results for holomorphic functions of several variables \cite{Ebel}. 
\begin{defn}
	Let 
	$\mathbb{R}_{>0}^l:=\{\mathbf{y}=(y_1, \dots, y_l)\in \mathbb{R}^l~|~ y_j>0 \text{~for all~}j \},$
	$\mathbf{r}=(r_1, \dots, r_l)\in \mathbb{R}^l_{>0}$, $\mathbf{a}\in \mathbb{C}^l$. Then, 
    $\Delta_{\mathbf{r}}(\mathbf{a}):=\{\mathbf{z}\in\mathbb{C}^l~|~|z_j-r_j|<r_j, 1\leq j\leq l \}$
    is called the polycylinder around $\mathbf{a}$ with (poly-)radius $\mathbf{r}$. The boundary of the closure of $\Delta_{\mathbf{r}}(\mathbf{a})$ contains an $n$-dimensional torus
    $T_{\mathbf{r}}(\mathbf{a}):=\{\mathbf{z}\in \mathbb{C}^l ~|~ |z_j-a_j|=r_j, 1\leq j\leq l \}.$
\end{defn}
In order to simplify the expressions in our proof, we introduce multiindices. Let $v_j$, $1\leq j\leq l$, be nonnegative integers, $\mathbf{z}=(z_1, \dots, z_l)\in\mathbb{C}^l$. Denote $\mathbf{v}:=(v_1, \dots, v_l)$, $|\mathbf{v}|=v_1+\cdots+v_l$,
$\mathbf{v}!:=v_1!\cdots v_l!$,  $\mathbf{z}^{\mathbf{v}}:=z_1^{v_1}\cdots z_l^{v_l}$, and
$$D^{\mathbf{v}}f=\frac{\partial^{|\mathbf{v}|}}{\partial z_1^{v_1}\cdots \partial z_l^{v_l}}.$$
We have the following result. 
\begin{lem}[\cite{Ebel}, Chapter 2, Propositions 2.7 and 2.11 ]\label{lem-Multi-Cauchy}
	Let $U\subset \mathbb{C}^l$ be open and $f: U\rightarrow\mathbb{C}$ holomorphic. Furthermore, let $\mathbf{w}\in U$ and $\Delta:=\Delta_{\mathbf{r}}(\mathbf{w})$ be a polycylinder around $\mathbf{w}$ with $\bar{\Delta}\subset U$, $T=T_{\mathbf{r}}(\mathbf{w})$. Then $f$ can be expanded as a power series 
	$$f(\mathbf{z})=\sum_{\mathbf{v}=0}^{\infty} a_{\mathbf{v}} (\mathbf{z}-\mathbf{w})^{\mathbf{v}}=\sum_{v_1\geq 0, \cdots, v_l\geq 0} a_{\mathbf{v}}(z_1-w_1)^{v_1}\cdots (z_l-w_l)^{v_l}$$
	in a neighborhood of $\mathbf{w}$, with coefficients
	\begin{equation*}
	a_{\mathbf{v}}=\frac{D^{\mathbf{v}}f(\mathbf{w})}{\mathbf{v}!}=\frac{1}{\mathbf{v}!}\frac{\partial^{v_1+\cdots+v_l}f}{\partial z_1^{v_1}\cdots \partial z_l^{v_l}}(\mathbf{w})=\(\frac{1}{2\pi i}\)^l \int_{T} \frac{f(\boldsymbol{\zeta})}{(\zeta_1-w_1)^{v_1+1}\cdots (\zeta_l-w_l)^{v_l+1}} d\boldsymbol{\zeta}.
	\end{equation*}
\end{lem}

\section{Proof of theorems}

\subsection{Associated Dirichlet series}

Let $(a, q)=1$. We define a function $\lambda_a(n)$ in the following way, 
\begin{equation}\label{defn-lambda-a}
\lambda_a(n)=\begin{cases}
1& \text{if~} n ~\text{square-free},~ p|n\Rightarrow p\equiv a \bmod q,\\
0 & \text{otherwise}.
\end{cases}
\end{equation}

We consider the Dirichlet series
\begin{equation}\label{F-a-defn}
F(s; a, z):=\sum_{n=1}^{\infty} \frac{ (z\lambda_a(n))^{\omega(n)}}{n^s}=\prod_p \(1+\frac{z\lambda_a(p)}{p^s} \), \quad (\Re(s)>1),
\end{equation}
where $\omega(n)$ is the number of distinct prime factors of $n$. Let $\chi_0$ be the principal character modulo $q$, denote
\begin{equation*}
F(s; z):=\sum_{n=1}^{\infty} \frac{\mu^2(n)(z \chi_0(n))^{\omega(n)}}{n^s}=\prod_{p} \( 1+\frac{z\chi_0(p)}{p^s}\)=\prod_{p\nmid q}\(1+\frac{z}{p^s} \), \quad (\Re(s)>1).
\end{equation*}
where $\mu(n)$ is the M\"{o}bius function.

Then we have the following lemma. 
\begin{lem}\label{lem-F-to-L}
	For any $A>0$, $|z|\leq A$, and $\Re(s)>1$, 
	\begin{equation*}
	F(s; a, z)=\(L(s, \chi_0)\)^{\frac{z}{\phi(q)}}\prod_{\chi\neq\chi_0} \(L(s, \chi)\)^{\frac{\bar{\chi}(a)z}{\phi(q)}} G_1(s; a, z),
	\end{equation*}
	and 
	\begin{equation*}
	F(s; z)=\( L(s, \chi_0)\)^{z} G_2(s; z),
	\end{equation*}
	where $\chi$ is a Dirichlet character modulo $q$, and $G_1(s; a, z)$ and $G_2(s; z)$ are absolutely convergent for $\Re(s)>\frac{1}{2}$. 
\end{lem}

\vspace{1em}

Given $\mathbf{a}=(a_1, a_2, \dots, a_k)$, assume there are $l$ distinct values $b_1, \dots, b_l$ in  the coordinates of $\mathbf{a}$. We assume $b_i$ ($1\leq i\leq l$) appears $k_i (>0)$ times in $\mathbf{a}$ with $k_1+k_2+\cdots +k_l=k$. 
Let $\mathbf{k}(\mathbf{a}):=(k_1, k_2, \dots, k_l)$, $\mathbf{b}(\mathbf{a}):=(b_1, \dots, b_l)$, and $\mathbf{z}=(z_1, z_2,\dots, z_l)$. Denote
\begin{equation}\label{F-general-defn}
F(s; \mathbf{a}, \mathbf{z}):=\prod_{j=1}^l F(s; b_j, z_j)=\prod_{j=1}^l \prod_p \(1+\frac{z_j\lambda_{b_j}(p)}{p^s} \).
\end{equation}
Let $\mathbf{n}=(n_1, \dots, n_l)\in \mathbb{Z}^l$ ($n_j>0, 1\leq j\leq l$). We write the Dirichlet series $F(s; \mathbf{a}, \mathbf{z})=\sum_{\mathbf{n>0}}\frac{a(\mathbf{n}; \mathbf{z})}{P^s(\mathbf{n})}$ with $P(\mathbf{n})=\prod_{1\leq j\leq l} n_j$. Then, 
\begin{equation*}
a(\mathbf{n}; \mathbf{z})=\sum_{\substack{\mathbf{k}=(k_1, \cdots, k_l)\\ k_j\geq 0}} c(\mathbf{k}, \mathbf{n}) z_1^{k_1}\cdots z_l^{k_l},
\end{equation*}
for some $c(\mathbf{k}, \mathbf{n})\in \mathbb{Z}^+$. Thus, for given $\mathbf{a}$, by Lemma \ref{lem-Multi-Cauchy}, 
\begin{equation}\label{Mk-integral}
M_k(x; \mathbf{a})=\sum_{P(\mathbf{n})\leq x} c(\mathbf{k}(\mathbf{a}), \mathbf{n})=\(\frac{1}{2\pi i}\)^l \oint_{|z_l|=r_l}\cdots\oint_{|z_1|=r_1} \(\sum_{P(\mathbf{n})\leq x} a(\mathbf{n}; \mathbf{z})\) \frac{dz_1}{z_1^{k_1+1}}\cdots\frac{dz_l}{z_l^{k_l+1}}.
\end{equation}

\subsection{A Uniform Result}
First, we prove the following result. 
\begin{theorem}\label{thm-general}
	For any $A>0$,  fixed $q\geq 3$ and fixed $l\geq 1$, uniformly for $2\leq k\leq A\log\log x$,  we have
	\begin{align*}
	M_k(x; \mathbf{a})=\frac{x}{\log x}\left\{\frac{1}{\phi(q)}Q_{\mathbf{k}}\(\frac{\log\log x}{\phi(q)}\) + O_{A, q, l}\(\frac{1}{\phi^k(q)}\frac{(\log\log x)^k}{k_1!\cdots k_l!\log x}\)\right\},
	\end{align*}
	where $Q_{\mathbf{k}}(X)$ is a polynomial of degree at most $k-1$ ($k=k_1+\cdots+k_l$).
	In particular, the coefficient of the term $\frac{x}{\log x}(\log\log x)^{k-1}$ is $\frac{1}{\phi^k(q)}\frac{k}{k_1!k_2!\cdots k_l!},$
	and the coefficient of $\frac{x}{\log x}(\log\log x)^{k-2}$ is  
	\begin{equation*}
	\frac{1}{\phi^k(q)}\frac{k(k-1)}{k_1!k_2!\cdots k_l!}\(\gamma+B+\frac{1}{k}\sum_{j=1}^k C(q, a_j)-\sum_{p|q}\frac{1}{p}\),
	\end{equation*}
	where $\gamma$ is Euler's constant, and $B:=\sum_p\(\log\(1-\frac{1}{p}\)+\frac{1}{p}\)$ is Mertens' constant. 
	
\end{theorem}

\noindent \textbf{\textit{Proof of Theorem \ref{thm-general}.}} By Lemma \ref{lem-F-to-L} and (\ref{F-general-defn}), we have
\begin{align}\label{F-gene-zeta-H}
F(s; \mathbf{a}, \mathbf{z})&=(L(s, \chi_0)^{\frac{z_1+\cdots+z_l}{\phi(q)}} \prod_{\chi\neq\chi_0}(L(s, \chi))^{\frac{\bar{\chi}(b_1)z_1+\cdots+\bar{\chi}(b_l)z_l}{\phi(q)}}\prod_{j=1}^l G_1(s; b_j, z_j)\nonumber\\
&=(\zeta(s))^{\frac{z_1+\cdots+z_l}{\phi(q)}}H(s; \mathbf{a}, \mathbf{z}),
\end{align}
where 
\begin{align*}
H(s; \mathbf{a}, \mathbf{z})&=\prod_{p|q}\(1-\frac{1}{p^s}\)^{\frac{z_1+\cdots+z_l}{\phi(q)}}\prod_{\chi\neq\chi_0}(L(s, \chi))^{\frac{\bar{\chi}(b_1)z_1+\cdots+\bar{\chi}(b_l)z_l}{\phi(q)}}\prod_{j=1}^l G_1(s; b_j, z_j)\nonumber\\
&=\prod_p \( 1-\frac{1}{p^s}\)^{\frac{z_1+\cdots+z_l}{\phi(q)}}\prod_{j=1}^l  \(1+\frac{z_j\lambda_{b_j}(p)}{p^s} \)
\end{align*}

Let $\sigma=\Re(s)$. Kolesnik \cite{Koles} showed that, for $\frac{1}{2}\leq \sigma\leq 1$,
\begin{equation}\label{L-estimate}
|L(s, \chi)|\ll (|t|+2)^{\frac{35}{108}(1-\sigma)}q^{1-\sigma} \log^3 (q(|t|+2)).
\end{equation}
Let $q$ be fixed. By Lemma \ref{lem-zerofree} and (\ref{L-estimate}), for any $A>0$, $|z_j|\leq A$ ($1\leq j\leq l$), and $0<\delta<1$, we can choose $c_0=c_0(A, \delta)$ such that, $L(s, \chi)$ has no zeros in the region $\sigma\geq 1-c_0/(\log (|t|+2))$, and by Theorem 11.4 in \cite{Mont-Vau}, in this region, 
$|H(s; \mathbf{a}, \mathbf{z})|\ll_{q, A, \delta} (|t|+2)^{1-\delta}.$
Thus, by Definitions \ref{defn-P} and \ref{defn-T}, $F(s; \mathbf{a}, \mathbf{z})$ is in  $\mathcal{T}(\frac{z_1+\cdots+z_l}{\phi(q)}, w; c_0, \delta, M)$. By (\ref{F-gene-zeta-H}) and following the same proof of Lemma \ref{lem-Dirichlet-z} (\cite{Tene}, Part II, Theorem 5.2, the only difference in the proof is the expansion of $H(s; \mathbf{a}, \mathbf{z})$), we deduce that, 
\begin{equation}\label{sum-a-z}
\sum_{n_1\cdots n_l\leq x} a(\mathbf{n}; \mathbf{z})=x(\log x)^{\frac{z_1+\cdots+z_l}{\phi(q)}-1}\left\{ u_0(\mathbf{a}; \mathbf{z})+O_A\(\frac{1}{\log x}\)\right\},
\end{equation}
where
\begin{equation}\label{pf-them-defn-u}
u_0(\mathbf{a}; \mathbf{z})=\frac{z_1+\cdots+z_l}{\phi(q)} u(\mathbf{a}; \mathbf{z}), \quad \text{with}\quad u(\mathbf{a}; \mathbf{z}):=\frac{H(1; \mathbf{a}, \mathbf{z})}{\Gamma\(\frac{z_1+\cdots+z_l}{\phi(q)}+1\)}.
\end{equation}

By (\ref{Mk-integral}), (\ref{sum-a-z}), and Lemma \ref{lem-Multi-Cauchy}, we have
\begin{align*}
M_k(x; \mathbf{a})=\frac{x}{\log x}\left\{\frac{1}{\phi(q)}Q_{ \mathbf{k}}\(\frac{\log\log x}{\phi(q)}\) + \widetilde{R}(x)\right\},
\end{align*}
where $Q_{\mathbf{k}}(X)$ is a polynomial of degree at most $k-1$ ($k=k_1+\cdots+k_l$), 
\begin{align}\label{pf-thm-general-Q0}
Q_{\mathbf{k}}(X)&:=\bigg\{\sum_{m_1+j_1=k_1-1}\sum_{m_2+j_2=k_2}\cdots\sum_{m_l+j_l=k_l}+\sum_{m_1+j_1=k_1}\sum_{m_2+j_2=k_2-1}\cdots\sum_{m_l+j_l=k_l}\nonumber\\
&\quad\quad +\cdots +\sum_{m_1+j_1=k_1}\cdots\sum_{m_{l-1}+j_{l-1}=k_{l-1}}\sum_{m_l+j_l=k_l-1}\bigg\}\nonumber\\
&\quad\quad \frac{1}{m_1!j_1!\cdots m_l!j_l!}\frac{\partial^{m_1+\cdots+m_l}}{\partial z_1^{m_1}\cdots \partial z_l^{m_l}}u(\mathbf{a}; (0, \cdots, 0)) X^{j_1+\cdots+j_l},
\end{align}
and 
\begin{equation}\label{pf-thm-R-est-eq1}
\widetilde{R}(x)\ll_A \frac{1}{(2\pi)^l \log x}\prod_{j=1}^l \oint_{|z_j|=r_j}  (\log x)^{\frac{\Re(z_j)}{\phi(q)}} \frac{|dz_j|}{|z_j|^{k_j+1}}.
\end{equation}
Taking $r_j=\frac{\phi(q)k_j}{\log\log x}$, we have
\begin{align}\label{pf-thm-R-est-eq2}
&\oint_{|z_j|=r_j}  (\log x)^{\frac{\Re(z_j)}{\phi(q)}} \frac{|dz_j|}{|z_j|^{k_j+1}}=\(\frac{\log\log x}{\phi(q)k_j} \)^{k_j} \int_0^{2\pi} e^{k_j\cos\theta}d\theta\nonumber\\
&\quad \leq \(  \frac{\log\log x}{\phi(q)k_j} \)^{k_j}\( 2\int_0^{\frac{\pi}{2}} e^{k_j\cos\theta}d\theta+\pi\)\nonumber\\
&\quad =\(  \frac{\log\log x}{\phi(q)k_j} \)^{k_j} \(2\int_0^1 e^{k_j t} \frac{dt}{\sqrt{1-t^2}}+\pi \)\nonumber\\
&\quad \leq \(  \frac{\log\log x}{\phi(q)k_j} \)^{k_j} \(2 e^{k_j} \int_0^1 e^{-k_j(1-t)}\frac{dt}{\sqrt{1-t}}+\pi \)\nonumber\\
&\quad \leq \(  \frac{\log\log x}{\phi(q)k_j} \)^{k_j}\( 2\Gamma\(\frac{1}{2}\)e^{k_j}k_j^{-\frac{1}{2}}+\pi\).
\end{align}
Substituting (\ref{pf-thm-R-est-eq2}) into (\ref{pf-thm-R-est-eq1}), we get
\begin{equation*}
\widetilde{R}_N(x)\ll_{A, l} \frac{1}{\phi^k(q)}\frac{(\log\log x)^k}{k_1!\cdots k_l!\log x}.
\end{equation*}

Theorem \ref{thm-general} follows. \qed
\vspace{1em}

\noindent \textbf{\textit{Remark 4.}} Similar to the proof of Lemma \ref{lem-Tene-Ck}, we write 
$$u(\mathbf{a}, \mathbf{z})=u(\mathbf{a}, \mathbf{r})+\sum_{|\mathbf{v}|=1}D^{\mathbf{v}}u(\mathbf{a}, \mathbf{r})+\sum_{|\boldsymbol{\beta}|=2}(\mathbf{z}-\mathbf{r})^{\boldsymbol{\beta}}R_{\boldsymbol{\beta}}(\mathbf{z}), $$ 
where
$$R_{\boldsymbol{\beta}}(\mathbf{z})=\frac{|\boldsymbol{\beta}|}{\boldsymbol{\beta}!}\int_0^1 (1-t)D^{\boldsymbol{\beta}}u(\mathbf{a}, \mathbf{r}+t(\mathbf{z}-\mathbf{r}))dt.$$
Then, by \eqref{Mk-integral} and \eqref{sum-a-z}, using a similar proof as in Lemma \ref{lem-Tene-Ck},  we have
\begin{equation}\label{pf-thm-general-uform}
M_k(x; \mathbf{a})=\frac{1}{\phi^k(q)}\frac{k}{k_1!k_2!\cdots k_l!}\frac{x(\log\log x)^{k-1}}{\log x}\left\{g\(\frac{\phi(q)}{\log\log x}; \mathbf{k} \)+O_{A, q, l}\(\frac{k}{(\log\log x)^2} \) \right\},
\end{equation}
where
$g(z; \mathbf{k}):=\sum_{j=1}^l \frac{k_j}{k}u(\mathbf{a}; (k_1 z, \cdots, k_{j-1} z, k'_j z, k_{j+1}z, \cdots, k_l z))$ with $k'_j=k_j-1$. Moreover, if $|kz|\leq A$, then $|g(z, \mathbf{k})|=O_{A, q, l}(1)$. 

\subsection{Proof of Theorems \ref{thm-same} and \ref{thm-general-quotient}}\label{sec-pf-thm-same}
For $\mathbf{a}=(a, \cdots, a)$, $(a, q)=1$, this is a special case of Theorem \ref{thm-general}. Denote
\begin{equation*}
H(s; a, z):=F(s; a, z)(\zeta(s))^{-\frac{z}{\phi(q)}}=\prod_p \(1-\frac{1}{p^s}\)^{\frac{z}{\phi(q)}}\(1+\frac{z\lambda_a(p)}{p^s} \), 
\end{equation*}
and
\begin{equation}\label{pf-thm-same-h-defn}
h(a; z):=\frac{H(1; a, z)}{\Gamma\(\frac{z}{\phi(q)}+1\)}. 
\end{equation}

Hence, for this special case, by \eqref{sum-a-z} and Lemma \ref{lem-Tene-Ck}, we get the following result. 

\begin{lem}\label{lem-Mk-same-pf}
For $\mathbf{a}=(a, \ldots, a)$ and any $A>0$, uniformly for $2\leq k\leq A\log\log x$, we have
\begin{align*}
M_k(x; \mathbf{a})=&\frac{1}{\phi^k(q)}\frac{x}{\log x} \frac{(\log\log x)^{k-1}}{(k-1)!}\bigg\{1+\frac{k-1}{\log\log x}C_{a, q}\nonumber\\
&+ \frac{(k-1)(k-2)\phi^2(q)}{(\log\log x)^2}\widetilde{h}\(a; \frac{(k-3)\phi(q)}{\log\log x}\)+O_{A, q}\(\frac{k^3}{(\log\log x)^4}\) \bigg\}, 
\end{align*}	
where
\begin{equation*}
C_{a, q}:=\phi(q)h'(a, 0)=\gamma+\sum_p \(\log\(1-\frac{1}{p}\)+\frac{\phi(q)\lambda_a(p)}{p} \),
\end{equation*}
$\gamma\approx 0.57722$ is Euler's constant, and 
$$\widetilde{h}(a, z)=\int_0^1 h''(a, tz)(1-t)dt.$$
\end{lem} 
\noindent \textbf{\textit{Remark 5.}} Notice that, for $|z|\leq A$, the function $|h''(a, z)|=O_{q, A}\( 1\)$ and $|h^{(4)}(a, z)|=O_{q, A}(1)$.

\bigskip

We also require a formula for $S_k(x)$. By Lemma \ref{lem-F-to-L}, and Definitions \ref{defn-P} and \ref{defn-T}, $F(s; z)$ is in $\mathcal{T}(z, w; c_0, \delta, M)$. Denote
\begin{equation*}
G(s; z):=F(s; z)(\zeta(s))^{-z}=\prod_p \(1-\frac{1}{p^s} \)^{z}\(1+\frac{z\chi_0(p)}{p^s}\), 
\end{equation*} 
and 
\begin{equation}\label{pf-thm-same-g-defn}
g(z):=\frac{G(1; z)}{\Gamma(z+1)}.
\end{equation}
Then, applying Lemma \ref{lem-Dirichlet-z} and Lemma \ref{lem-Tene-Ck} successively, we get the following lemma. 
\begin{lem}\label{lem-Sk-pf}
	For any $A>0$, uniformly for $2\leq k\leq A\log\log x$, we have
\begin{align*}
S_k(x)=&\frac{x}{\log x}\frac{(\log\log x)^k}{(k-1)!}\bigg\{1+\frac{k-1}{\log\log x}g'\(0\)+\frac{(k-1)(k-2)}{(\log\log x)^2}\widetilde{g}\(\frac{k-3}{\log\log x}\)\nonumber\\
&+O_{A, q}\(\frac{k^3}{(\log\log x)^4}\) \bigg\},
\end{align*}
where $g'(0)=\gamma+B-\sum_{p|q}\frac{1}{p}$, $\gamma$ is Euler's constant,  $B=\sum_p \( \log\(1-\frac{1}{p}\)+\frac{1}{p}\)$ is Mertens' constant in (\ref{Mertens-thm}), and 
$$\widetilde{g}(z)=\int_0^1 g''(tz)(1-t)dt.$$
\end{lem}

\noindent \textbf{\textit{Remark 6.}} Here for $|z|\leq A$, the function $|g''(z)|=O_{q, A}(1)$ and $|g^{(4)}(z)|=O_{q, A}(1)$.
\vspace{1em}

\noindent \textbf{\textit{Proof of Theorem \ref{thm-same}.}}
By Lemmas \ref{lem-Mk-same-pf} and \ref{lem-Sk-pf}, we get
\begin{align*}
M_k(x, \mathbf{a})-\frac{1}{\phi^k(q)}S_k(x)&=\frac{1}{\phi^k(q)}\frac{x}{\log x}\frac{(\log\log x)^{k-2}}{(k-2)!}\bigg\{C(q, a) +\frac{k-2}{\log\log x}\phi^2(q)\widetilde{h}\(a; \frac{(k-3)\phi(q)}{\log\log x}\)\nonumber\\
&-\frac{k-2}{\log\log x}\widetilde{g}\(\frac{k-3}{\log\log x}\) +O_{A, q}\( \frac{k}{(\log\log x)^2} \) \bigg\}.
\end{align*}
For the cases of fixed $k$ and $k=o(\log\log x)$, by Remarks 5 and 6, and Lemma \ref{lem-Sk-pf}, we immediately get the conclusions in Theorem \ref{thm-same} using the equality
\begin{equation*}
\frac{M_k(x, \mathbf{a})}{\frac{1}{\phi^k(q)}S_k(x)}=1+\frac{M_k(x, \mathbf{a})-\frac{1}{\phi^k(q)}S_k(x)}{\frac{1}{\phi^k(q)}S_k(x)}.
\end{equation*}

For any fixed $A>0$, if $k\sim A\log\log x$, by Lemmas \ref{lem-Mk-same-pf} and \ref{lem-Sk-pf}, and (\ref{pf-thm-general-uform}), as $x\rightarrow \infty$, the above quotient will approach $$\frac{h(a, A\phi(q))}{g(A)}=\prod_p \frac{1+\frac{A\phi(q)\boldsymbol{\mathbbm{1}}_{p\equiv a \bmod q}(p)}{p}}{1+\frac{A\chi_0(p)}{p}},$$ 
where $h(a, z)$ and $g(z)$ are defined in (\ref{pf-thm-same-h-defn}) and (\ref{pf-thm-same-g-defn}), respectively.  \qed
\vspace{1em}

\noindent \textbf{\textit{Proof of Theorem \ref{thm-general-quotient}.}}
For fixed $k$, by Theorem \ref{thm-general} and Lemma \ref{lem-Sk-pf}, we have
\begin{align*}
&M_k(x; \mathbf{a})-\frac{1}{\phi^k(q)}\frac{k!}{k_1!k_2!\cdots k_l!} S_k(x)\nonumber\\
&\qquad =\frac{1}{\phi^k(q)}\frac{k(k-1)}{k_1!k_2!\cdots k_l!}\frac{x}{\log x}(\log\log x)^{k-2}\left\{\frac{1}{k}\sum_{j=1}^k C(q, a_j)+O_{k, q, l}\(\frac{1}{\log\log x} \)  \right\}.
\end{align*}
Thus, 
\begin{align*}
\frac{M_k(x; \mathbf{a})}{\frac{1}{\phi^k(q)}\frac{k!}{k_1!k_2!\cdots k_l!} S_k(x)}&=1+\frac{M_k(x; \mathbf{a})-\frac{1}{\phi^k(q)}\frac{k!}{k_1!k_2!\cdots k_l!} S_k(x)}{\frac{1}{\phi^k(q)}\frac{k!}{k_1!k_2!\cdots k_l!} S_k(x)}\\
&= 1+\frac{k-1}{\log\log x} \frac{1}{k}\sum_{j=1}^k C(q, a_j)+O_{q, k, l}\(\frac{1}{(\log\log x)^2} \).
\end{align*}

For any fixed $A>0$, if $k\sim A\log\log x$  and $e_j:=\lim_{x\rightarrow\infty} \frac{k_j}{\log\log x}$ exists, by (\ref{pf-thm-general-uform}) and Lemma \ref{lem-Sk-pf}, as $x\rightarrow\infty$, the above quotient will approach 
$$\frac{u(\mathbf{a}; ( \phi(q)e_1, \cdots,  \phi(q)e_l))}{g( A)}=\prod_p \frac{\prod_{j=1}^l \(1+\frac{\phi(q)e_j \boldsymbol{\mathbbm{1}}_{p\equiv b_j \bmod q}(p)}{p}\)}{1+\frac{A\chi_0(p)}{p}},$$
where $u(\mathbf{a}; \mathbf{z})$ and $g( z)$ are defined in (\ref{pf-them-defn-u}) and (\ref{pf-thm-same-g-defn}) respectively.  \qed

\vspace{1em}
\textbf{Acknowledgement.} This research is partially
supported by NSF grant DMS-1501982. I would like to thank my advisor, Professor Kevin Ford, for his useful comments and financial support to finish this project. I am grateful for the helpful comments and advice of Prof. Andrew Granville. I'd like to thank Dr. Nathan McNew and Prof. Carl Pomerance for the helpful discussions on the least prime in arithmetic progressions. The author thanks Peter Humphries for his comments and for pointing out the reference \cite{Norton}. The author would like to thank the referee for his/her careful reading and helpful comments, and for pointing out the reference \cite{Baker}.

{\footnotesize
Department of Mathematics, University of Illinois at Urbana-Champaign, 1409 West Green Street, Urbana, IL 61801, USA

\textit{E-mail}:  xmeng13@illinois.edu,

\end{document}